\documentclass[11pt,english]{article}
\usepackage{palatino}
\usepackage[T1]{fontenc}
\usepackage[latin1]{inputenc}
\usepackage{a4}
\usepackage{amssymb}

\makeatletter

\providecommand{\boldsymbol}[1]{\mbox{\boldmath $#1$}}

 \usepackage{amsthm, amsmath, amsfonts, amssymb}
 \theoremstyle{plain}    
 \newtheorem{thm}{Theorem}
 \theoremstyle{plain}    
 \newtheorem{lem}[thm]{Lemma} 
 \theoremstyle{plain}    
 \newtheorem{prop}[thm]{Proposition} 

 \theoremstyle{definition}
 \newtheorem{defn}[thm]{Definition}

\usepackage{babel}
\makeatother
\begin{document}

\title{Regular deformations of completely integrable systems}

\author{Nicolas Roy%
\thanks{\noindent Address : Geometric Analysis Group, Institut für Mathematik,
Humboldt Universität, Rudower Chaussee 25, Berlin D-12489, Germany.
Email : \texttt{roy@math.hu-berlin.de}%
}~%
\thanks{The author would like to thank Janko Latschev for reading the manuscript
and the referee for his very useful critics.%
}~%
\thanks{This paper has been partially supported by the European Commission
through the Research Training Network HPRN-CT-1999-00118  \char`\"{}Geometric
Analyis\char`\"{}.%
}}

\maketitle
\begin{abstract}
We study several aspects of the regular deformations of completely
integrable systems. Namely, we prove the existence of a Hamiltonian
normal form for these deformations and we show the necessary and sufficient
conditions a perturbation has to satisfy in order for the perturbed
Hamiltonian to be a first order deformation.
\end{abstract}

\section*{Introduction}

This article presents some results concerning the deformations of
regular completely integrable (CI in short) systems. These are the
dynamical systems defined by a Hamiltonian $H_{0}\in C^{\infty}\left(\mathcal{M}\right)$
on a symplectic manifold $\mathcal{M}$ admitting a momentum map,
i.e. a collection $\boldsymbol{A}=\left(A_{1},...,A_{d}\right):\mathcal{M}\rightarrow\mathbb{R}^{d}$
of $d$ smooth functions, $d$ being half of the dimension of $\mathcal{M}$,
satisfying $\left\{ A_{j},H\right\} =0$ and $\left\{ A_{j},A_{k}\right\} =0$
for all $j,k:1...d$, and whose differentials $dA_{j}$ are linearly
independent almost everywhere. Then, the Arnol'd-Mineur-Liouville
Theorem \cite{arnold_1,mineur,liouville} insures that in a neighborhood
of any connected component of any compact regular fiber $\boldsymbol{A}^{-1}\left(a\right)$,
$a\in\mathbb{R}^{d}$, of the momentum map, there exists a fibration
in Lagrangian tori along which $H_{0}$ is constant. These tori are
thus invariant by the dynamics generated by the associated Hamiltonian
vector field $X_{H}$. 

Despite the {}``local'' character of the Arnol'd-Mineur-Liouville
Theorem, it is tempting to try to glue together these {}``local''
fibrations in the case of \emph{regular} CI Hamiltonians, i.e. those
for which there exists, near each point of \emph{}$\mathcal{M}$,
a local fibration in invariant Lagrangian tori. Unfortunately, this
is not always possible. Some CI Hamiltonians do not admit any (global)
fibration in Lagrangian tori and some others admit several different
ones%
\footnote{For example, the free particle moving on the sphere $S^{2}$ is such
a system.%
}. Nevertheless, these examples belong to the non-generic (within the
class of regular CI Hamiltonians) class of \emph{degenerate} Hamiltonians
and one can show that imposing a nondegeneracy condition insures that
there exists a fibration of $\mathcal{M}$ in Lagrangian tori along
which $H_{0}$ is constant, and moreover that it is unique. The genericity
of nondegeneracy conditions thus motivates the study of fibrations
in Lagrangian tori $\mathcal{M}\overset{\pi}{\rightarrow}\mathcal{B}$.
Such a fibration actually gives rise to several natural geometric
structures which we review in the first section. 

Starting from a regular CI Hamiltonian $H_{0}\in C^{\infty}\left(\mathcal{M}\right)$,
it is well-known since Poincaré's work \cite{poincare_123} that adding
a small perturbation $\varepsilon H_{1}$ will destroy its integrable
character and yield chaotic behaviours. Nevertheless, it is relevant
to investigate the space of all CI Hamiltonians, since they are the
starting point of any perturbation theory, like the celebrated K.A.M.
Theory \cite{kolmogorov,arnold_2,moser_1} which tells us that one
can actually say a lot about the perturbed Hamiltonian $H_{\varepsilon}=H_{0}+\varepsilon H_{1}$
when $\varepsilon$ is small.

A first step towards the understanding of the space of all CI systems,
is to restrict ourselves to regular deformations of regular CI hamiltonians,
i.e smooth families of Hamiltonians $H_{\varepsilon}$ which are CI
and regular for each $\varepsilon$. After introducing a few necessary
tools in Section \ref{sec_ND}, we prove in Section \ref{sec_normal_form}
a normal form for regular deformations of CI Hamiltonians. Finally,
Section \ref{sec_first_order_deformations} is devoted to the study
of the first order deformations. We give there the condition on the
perturbation $H_{1}$ for the perturbed system $H_{0}+\varepsilon H_{1}$
to be CI up to $\varepsilon^{2}$.

\section{Geometric structures of regular CI systems}

In this section, we review several geometric structures which are
naturally associated with any fibration in Lagrangian tori $\mathcal{M}\overset{\pi}{\rightarrow}\mathcal{B}$.
In particular, we show that there exists a natural process of averaging
any tensor field in the direction of the fibers. This process then
allows us to prove (Proposition \ref{prop_decomposition_vectors})
that each symplectic vector field splits into two parts : the first
is Hamiltonian and the second is symplectic and preserves the fibration.
This will be used in Section \ref{sec_normal_form} to prove the Hamiltonian
normal form (Theorem \ref{theo_deformation_hamiltonian}) for regular
deformations.

First, let us fix some basic notations. We denote by $\mathcal{V}\left(\mathcal{M}\right)$
the space of smooth vector fields on the manifold $\mathcal{M}$.
A symplectic form $\omega$ on $\mathcal{M}$ provides a isomorphism
$\omega:\mathcal{V}\left(\mathcal{M}\right)\rightarrow\Omega^{1}\left(\mathcal{M}\right)$,
also denoted by $\omega$, i.e. $\omega\left(X\right)=\omega\left(X,.\right)$
for each $X\in C^{\infty}\left(\mathcal{M}\right)$. The inverse is
denoted by $\omega^{-1}:\Omega^{1}\left(\mathcal{M}\right)\rightarrow\mathcal{V}\left(\mathcal{M}\right)$.
For each vector field, we denote by $\phi_{X}^{t}$ its flow at time
$t$. Let $O\subset\mathcal{M}$ be any subset. We say that a vector
field $X$ is symplectic (resp. Hamiltonian) in $\mathcal{O}$ if
its associated $1$-form $\omega\left(X\right)$ is closed (resp.
exact) in $\mathcal{O}$. To each Hamiltonian $H\in C^{\infty}\left(\mathcal{M}\right)$
we can associate a vector field $X_{H}=-\omega^{-1}\left(dH\right)$.
Now, given a fibration $\mathcal{M}\overset{\pi}{\rightarrow}\mathcal{B}$,
we say that a vector field $\tilde{X}\in\mathcal{V}\left(M\right)$
is a lift of a vector field $X\in\mathcal{V}\left(\mathcal{B}\right)$
if for each $b\in\mathcal{B}$ and each $m\in\pi^{-1}\left(b\right)$
we have $\pi_{*}\left(\widetilde{X}_{m}\right)=X_{b}$.

\subsection{The period bundle\label{sec_period_bundle}}

Let $\left(\mathcal{M},\omega\right)$ be a symplectic manifold of
dimension $2d$ and $\mathcal{M}\overset{\pi}{\rightarrow}\mathcal{B}$
a locally trivial fibration in Lagrangian tori, whose fibers are denoted
by $\mathcal{M}_{b}=\pi^{-1}\left(b\right)$, $b\in\mathcal{B}$.
The tangent spaces $L_{m}=T_{m}\mathcal{M}_{\pi\left(m\right)}$ of
the fibers form an integrable vector subbundle $L=\bigcup_{m\in\mathcal{M}}L_{m}$
of $T\mathcal{M}$. A theorem due to Weinstein \cite{weinstein1}
insures that each leaf of a Lagrangian foliation (not necessarily
a fibration) is naturally endowed with an affine structure. This affine
structure on a leaf $\mathcal{M}_{b}$ can actually be expressed in
a very convenient way (see \cite{woodhouse1}) in terms of the torsion-free
and flat connection $\nabla:\mathcal{V}\left(\mathcal{M}_{b}\right)\times\mathcal{V}\left(\mathcal{M}_{b}\right)\rightarrow\mathcal{V}\left(\mathcal{M}_{b}\right)$
defined by \[
\nabla_{X}Y=\omega^{-1}\left(\widetilde{X}\lrcorner d\left(\widetilde{Y}\lrcorner\omega\right)\right),\]
where $\widetilde{X}\in\Gamma\left(L\right)$ and $\widetilde{Y}\in\Gamma\left(L\right)$
extend $X$ and $Y$ in $\mathcal{V}\left(\mathcal{M}\right)$ and
are everywhere tangent to $L$. We denote by $\mathcal{V}_{\nabla}\left(\mathcal{M}_{b}\right)$
the space of parallel vector fields on $\mathcal{M}_{b}$. One can
see easily from the definition of $\nabla$ that a vector field $X\in\mathcal{V}\left(\mathcal{M}\right)$
is vertical and parallel on each fiber if and only if its associated
$1$-form $\omega\left(X\right)$ is a pull-back of a $1$-form on
$\mathcal{B}$.

Now, since the foliation under consideration actually defines a fibration,
the holonomy of $\nabla$ must vanish. Indeed, for each $b\in\mathcal{B}$,
any collection of smooth functions $f_{1},...,f_{d}\in C^{\infty}\left(\mathcal{B}\right)$
whose differentials $df_{j}$ are linearly independent near $b$,
provides $d$ Hamiltonian vector fields $X_{f_{1}\circ\pi},...,X_{f_{d}\circ\pi}\in\mathcal{V}\left(\mathcal{M}\right)$
everywhere tangent to the fibers, parallel on each fiber and linearly
independent in a neighborhood of $\mathcal{M}_{b}$. Therefore, they
form a global parallel frame on $\mathcal{M}_{b}$, implying that
the holonomy of $\nabla$ vanishes and that each fiber $\mathcal{M}_{b}$
is endowed with the structure of a standard%
\footnote{Here, {}``standard'' means holonomy-free. %
} affine torus. This implies that the space $\mathcal{V}_{\nabla}\left(\mathcal{M}_{b}\right)$
is a $d$-dimensional vector space and that the union $\bigcup_{b\in\mathcal{B}}\mathcal{V}_{\nabla}\left(\mathcal{M}_{b}\right)$
is naturally endowed with a structure of a smooth vector bundle over
$\mathcal{B}$. 

Since each fiber $\mathcal{M}_{b}$ is isomorphic to the standard
torus $\mathbb{T}^{d}$, we can consider among the parallel vector
fields on $\mathcal{M}_{b}$, those whose dynamics is $1$-periodic.
We denote this set by $\Lambda_{b}$. It is then easy to prove that
it is a lattice in $\mathcal{V}_{\nabla}\left(\mathcal{M}_{b}\right)$.
We call it the \emph{period lattice}. The genuine geometric content
of the Arnol'd-Mineur-Liouville Theorem \cite{arnold_1,mineur,liouville},
which is often hidden by the formulation in coordinates, amounts to
saying that the union $\Lambda=\bigcup_{b\in\mathcal{B}}\Lambda_{b}$,
called the \emph{period bundle,} is a smooth lattice subbundle of
$\bigcup_{b\in\mathcal{B}}\mathcal{V}_{\nabla}\left(\mathcal{M}_{b}\right)$.
This can be proved by constructing explicit smooth sections of this
bundle which are $1$-periodic, namely Hamiltonian vector fields $X_{\xi\circ\pi}$
where the function $\xi\in C^{\infty}\left(\mathcal{B}\right)$ is
called \emph{action} and given by $b\rightarrow\xi\left(b\right)=\int_{\gamma\left(b\right)}\theta$,
with $\theta$ any symplectic potential and $b\rightarrow\gamma\left(b\right)$
a smooth family of vertical cycles. Furthermore, this shows that smooth
(local) sections of $\Lambda$ are Hamiltonian.

The smoothness of the period bundle $\Lambda$ provides a way to relate
the spaces $\mathcal{V}_{\nabla}\left(\mathcal{M}_{b}\right)$ for
neighboring points $b$ and implies the existence of a natural integer
flat connection on the vector bundle $\bigcup_{b\in\mathcal{B}}\mathcal{V}_{\nabla}\left(\mathcal{M}_{b}\right)$.
This connection may have a non-vanishing holonomy, called the \emph{monodromy}.
Now, the symplectic form $\omega$ provides an isomorphism between
the sections of $\bigcup_{b\in\mathcal{B}}\mathcal{V}_{\nabla}\left(\mathcal{M}_{b}\right)$
and those of $T^{*}\mathcal{B}$. This gives the base space $\mathcal{B}$
a natural structure of an affine space, as discovered by Duistermaat
in \cite{duistermaat}.

\subsection{The torus action bundle}

Our discussion so far shows that given a fibration in Lagrangian tori
$\mathcal{M}\overset{\pi}{\rightarrow}\mathcal{B}$, there exists
a natural associated torus bundle acting on it. Indeed, for each $b\in\mathcal{B}$,
the quotient \[
\mathcal{G}_{b}=\mathcal{V}_{\nabla}\left(\mathcal{M}_{b}\right)/\Lambda_{b}\]
 is a Lie group isomorphic to the torus $\mathbb{T}^{d}$. This isomorphism
is not canonical, but it can be realised by choosing a basis of $\Lambda_{b}$.
We will denote the elements of $\mathcal{G}_{b}$ by $\left[X_{b}\right]$,
with $X_{b}\in\mathcal{V}_{\nabla}\left(\mathcal{M}_{b}\right)$,
since they are equivalence classes. Taking the union over all $b$,
we get a torus bundle $\mathcal{G}=\bigcup_{b\in\mathcal{B}}\mathcal{G}_{b}$.
It is a smooth bundle since the period bundle $\Lambda$ is so. We
stress the fact that $\mathcal{G}$ is in general \emph{not} a principal
bundle since there might not exist any global action of $\mathbb{T}^{d}$
on $\mathcal{G}$, because of the presence of monodromy which precisely
prevents us from choosing a global basis of $\Lambda$. On the other
hand, there exists a distinguished global section, since each fiber
is a group with a well-defined identity element. 

Although we cannot apply the general theory of connections on principal
bundles, there is a natural way to speak about local parallel sections
of $\mathcal{G}$ over a subset $\mathcal{O}\subset\mathcal{B}$.
These sections are simply local sections $b\rightarrow\left[X_{b}\right]$
of $\mathcal{G}$, with $b\rightarrow X_{b}$ being a local parallel
section of $\bigcup_{b\in\mathcal{B}}\mathcal{V}_{\nabla}\left(\mathcal{M}_{b}\right)$.
We denote the set of local parallel sections by $\Gamma_{\nabla}\left(\mathcal{O},\mathcal{G}\right)$.

\begin{lem}
\label{lem_local_section_torus}For each simply connected subset $\mathcal{O}\subset\mathcal{B}$,
the space $\Gamma_{\nabla}\left(\mathcal{O},\mathcal{G}\right)$ is
a Lie group isomorphic to the torus $\mathbb{T}^{d}$.
\end{lem}
\begin{proof}
If $\mathcal{O}$ is simply connected, then the monodromy vanishes
in $\mathcal{O}$ and there exist local sections $X_{1},...,X_{d}\in\Gamma\left(\mathcal{O},\Lambda\right)$
with $\left\{ X_{j}\left(b\right)\right\} $ generating the lattice
$\Lambda_{b}$ at each $b\in\mathcal{O}$. To each element $\left(t_{1},...t_{d}\right)\in\mathbb{T}^{d}=\mathbb{R}^{d}/\mathbb{Z}^{d}$,
we associate $\left[X\right]=\left[t_{1}X_{1}+...+t_{d}X_{d}\right]\in\Gamma_{\nabla}\left(\mathcal{O},\mathcal{G}\right)$.
One easily verifies that this provides an isomorphism.
\end{proof}
Let us now describe how the bundle $\mathcal{G}$ acts on $\mathcal{M}$.
First, for each $b$ the group $\mathcal{G}_{b}$ acts naturally on
$\mathcal{M}_{b}$ in the following way. \begin{eqnarray*}
\mathcal{G}_{b}\times\mathcal{M}_{b} & \rightarrow & \mathcal{M}_{b}\\
\left(\left[X_{b}\right],m\right) & \rightarrow & \left[X_{b}\right]\left(m\right)=\phi_{X_{b}}^{1}\left(m\right),\end{eqnarray*}
where $X_{b}\in\mathcal{V}_{\nabla}\left(\mathcal{M}_{b}\right)$
is a representative of the class $\left[X_{b}\right]$. One can see
easily that this action is commutative, free, transitive and affine
with respect to Weinstein's connection on $\mathcal{M}_{b}$. Now,
given any section $g\in\Gamma\left(\mathcal{G}\right)$, its restriction
$\left.g\right|_{\mathcal{O}}$ to any simply connected subset $\mathcal{O}\subset\mathcal{B}$
is of the form $\left.g\right|_{\mathcal{O}}=\left[X\right]$, where
$X\in\Gamma\left(\mathcal{O},\bigcup_{b\in\mathcal{B}}\mathcal{V}_{\nabla}\left(\mathcal{M}_{b}\right)\right)$.
We can then extend the previous fiberwise action of the groups $\mathcal{G}_{b}$
to a vertical action of the sections of the toric bundle $\mathcal{G}$
on $\mathcal{M}$ by \begin{eqnarray*}
\Gamma\left(\mathcal{G}\right)\times\mathcal{M} & \rightarrow & \mathcal{M}\\
\left(g,m\right) & \rightarrow & \left[X\right]\left(m\right)=\phi_{X}^{1}\left(m\right),\end{eqnarray*}
where $X\in\Gamma\left(\mathcal{O},\bigcup_{b\in\mathcal{B}}\mathcal{V}_{\nabla}\left(\mathcal{M}_{b}\right)\right)$
for any simply connected neighborhood $\mathcal{O}$ of $b=\pi\left(m\right)$.
This is well-defined since another representative $X^{'}$ of the
class of $\left[X\right]$ would differ from $X$ only by an element
of $\Gamma\left(\mathcal{O},\Lambda\right)$ which would provide $\phi_{X^{'}-X}^{1}=\mathbb{I}$.
This action naturally inherits the properties of the fiberwise action
and we can show that the following additional property holds when
we restrict ourselves to the parallel sections of $\mathcal{G}$.

\begin{lem}
For any simply connected subset $\mathcal{O}\subset\mathcal{B}$,
the group $\Gamma_{\nabla}\left(\mathcal{O},\mathcal{G}\right)$ acts
vertically on $\mathcal{M}$ in a symplectic way. 
\end{lem}
We call this action the \emph{toric action} of $\mathcal{G}$ on $\mathcal{M}$.
Even if this action is local, it provides a way to average any tensor
field on \emph{$\mathcal{M}$}. Indeed, according to Lemma \ref{lem_local_section_torus},
$\Gamma_{\nabla}\left(\mathcal{O},\mathcal{G}\right)$ is a compact
Lie group provided $\mathcal{O}\subset\mathcal{B}$ is simply connected.
It is thus endowed with its Haar measure $\mu_{\mathcal{G}}$ and
for any tensor field $T$ of any type on $\mathcal{M}$, we can define
its \emph{vertical average} $\left\langle T\right\rangle $ in the
following way. For each $m\in\mathcal{M}$, we set \[
\left\langle T\right\rangle _{m}=\int_{\Gamma_{\nabla}\left(\mathcal{O},\mathcal{G}\right)}\left(\left(\phi_{X}^{1}\right)_{*}T\right)_{m}d\mu_{\mathcal{G}},\]
where $\mathcal{O}\subset\mathcal{B}$ is any simply connected neighborhood
of $b=\pi\left(m\right)$. We can check that this definition does
not depend on the choice of $\mathcal{O}$. Choosing a basis $X_{1},...,X_{d}$
of $\Gamma\left(\mathcal{O},\Lambda\right)$ provides an explicit
expression for the averaged tensor, namely \[
\left\langle T\right\rangle _{m}=\int_{0}^{1}dt_{1}...\int_{0}^{1}dt_{d}\left(\left(\phi_{X_{1}}^{t_{1}}\right)_{*}\circ...\circ\left(\phi_{X_{d}}^{t_{d}}\right)_{*}T\right)_{m}.\]

A tensor field $T$ is called \emph{invariant under the toric action
of} $\mathcal{G}$, or simply $\mathcal{G}$\emph{-invariant}, if
for each local parallel section $X\in\Gamma_{\nabla}\left(\mathcal{O},\mathcal{G}\right)$
we have $\left(\phi_{X}^{1}\right)_{*}\left(T\right)=T$, or equivalently
$\mathcal{L}_{X}T=0.$ The following properties can be proved in a
straightforward way.

\begin{lem}
\label{lem_moyenne_action_torique_2}We have the following basic properties
:
\begin{enumerate}
\item $T$ is $\mathcal{G}$-invariant if and only if $\left\langle T\right\rangle =T$.
\item $\left\langle \left\langle T\right\rangle \right\rangle =\left\langle T\right\rangle $.
\item Each $p$-form $\alpha\in\Omega^{p}\left(\mathcal{M}\right)$ verifies
$\left\langle d\alpha\right\rangle =d\left\langle \alpha\right\rangle $.
\item Let $T$ and $S$ be two tensor fields. If $T$ is $\mathcal{G}$-invariant,
then the contraction $T\lrcorner S$ with respect to any two indices
verifies $\left\langle T\lrcorner S\right\rangle =T\lrcorner\left\langle S\right\rangle $. 
\item In particular, if $X\in\mathcal{V}\left(\mathcal{M}\right)$ is a
vector field and $\alpha=\omega\left(X\right)$ its associated $1$-form,
then we have $\omega\left(\left\langle \alpha\right\rangle \right)=\left\langle X\right\rangle $.
\end{enumerate}
\end{lem}

\subsection{Decomposition of symplectic vector fields}

The averaging process presented in the previous section provides a
way to decompose any symplectic vector field into the sum of a Hamiltonian
vector field and a symplectic vector field preserving the fibration.
The key step is the following lemma.

\begin{lem}
\label{lem_forme_ferme_moyenne_null_exacte}If $\alpha$ is a closed
$1$-form on $\mathcal{M}$ whose vertical average vanishes, then
it is exact. Moreover, one can choose the primitive $f\in C^{\infty}\left(\mathcal{M}\right)$,
$\alpha=df$, with the property $\left\langle f\right\rangle =0$.
\end{lem}
\begin{proof}
Let us work locally in a simply connected subset $\mathcal{O}\subset\mathcal{B}$.
There exists a basis $\left(X_{1},\cdots,X_{d}\right)$ of $\Gamma\left(\mathcal{O},\Lambda\right)$.
Choosing an {}``initial point'' $m\left(b\right)$ depending smoothly
on $b\in\mathcal{O}$, i.e. a smooth section of the restricted bundle
$\pi^{-1}\left(\mathcal{O}\right)\overset{\pi}{\rightarrow}\mathcal{O}$,
let us consider the smooth family of cycles $\gamma_{j}\left(b\right)$
consisting of the orbits $t\rightarrow\phi_{X_{j}}^{t}\left(m\left(b\right)\right)$.
The homology classes $\left[\gamma_{j}\left(b\right)\right]$ form
for each $b\in\mathcal{O}$ a basis of $H_{1}\left(\mathcal{M}_{b}\right)$.
On the other hand, since the fibration $\mathcal{M}\overset{\pi}{\rightarrow}\mathcal{B}$
is locally trivial and $\mathcal{O}$ is simply connected, the classes
$\left[\gamma_{j}\left(b\right)\right]$ form a basis of the homology
of $\tilde{\mathcal{O}}=\pi^{-1}\left(\mathcal{O}\right)$. 

Then, we show that for each $j=1..d$ and each $b\in\mathcal{O}$,
one has $\int_{\gamma_{j}\left(b\right)}\left\langle \alpha\right\rangle =\int_{\gamma_{j}\left(b\right)}\alpha$.
Indeed, one has \[
\int_{\gamma_{j}\left(b\right)}\left\langle \alpha\right\rangle =\int_{0}^{1}dt\,\left\langle \alpha\right\rangle \left(X_{j}\right)\circ\phi_{X_{j}}^{t}\left(m\left(b\right)\right)=\int_{0}^{1}dt\, X_{j}\lrcorner\left(\phi_{X_{j}}^{-t}\right)_{*}\left\langle \alpha\right\rangle .\]
Moreover, expressing the average $\left\langle \alpha\right\rangle $
in terms of the generators $X_{j}$, one obtains \[
\int_{\gamma_{j}\left(b\right)}\left\langle \alpha\right\rangle =\int_{0}^{1}dt_{1}...\int_{0}^{1}dt_{d}\int_{0}^{1}dt\,\left(\phi_{X_{1}}^{t_{1}}\right)_{*}\circ\cdots\circ\widehat{\left(\phi_{X_{j}}^{t_{j}}\right)}_{*}\circ\cdots\circ\left(\phi_{X_{d}}^{t_{d}}\right)_{*}\left(X_{j}\lrcorner\left(\phi_{X_{j}}^{t_{j}-t}\right)_{*}\alpha\right),\]
where the entry below $\widehat{\,\,}$ has been omitted. Then, we
check with a trivial change of variable that \[
\int_{0}^{1}dt_{j}\int_{0}^{1}dt\left(X_{j}\lrcorner\left(\phi_{X_{j}}^{t_{j}-t}\right)_{*}\alpha\right)=\int_{\gamma_{j}\left(b\right)}\alpha.\]
This implies that $\int_{\gamma_{j}\left(b\right)}\left\langle \alpha\right\rangle =\int_{\gamma_{j}\left(b\right)}\alpha$.

Finally, the hypothesis $\left\langle \alpha\right\rangle =0$ yields
$\int_{\gamma_{j}\left(b\right)}\alpha=0$, where the classes $\left[\gamma_{j}\left(b\right)\right]$
form a basis of the homology of $\tilde{\mathcal{O}}=\pi^{-1}\left(\mathcal{O}\right)$,
as shown before. Since $\alpha$ is closed, this implies that it is
actually exact. Thus, there exists a function $f\in C^{\infty}\left(\tilde{\mathcal{O}}\right)$
such that $\alpha=df$ in $\mathcal{\tilde{O}}$. This function is
unique up to a constant. On the other hand, we deduce from the property
$\left\langle df\right\rangle =d\left\langle f\right\rangle $ and
the hypothesis $\left\langle \alpha\right\rangle =0$ that $\left\langle f\right\rangle $
is a constant function. This allows us to choose the primitive $f$
in a unique way by requiring that $\left\langle f\right\rangle =0$.
This criterion is independent of the choice of the basis $\left(X_{1},...,X_{d}\right)$
and thus allows us to find a primitive $f$ of $\alpha$ globally
defined on $\mathcal{M}$.
\end{proof}
We need also the following property, which will be proved later in
the slightly more general case of time-dependent vector fields (Lemma
\ref{lem_vecteur_symplec_et_lift}). 

\begin{lem}
\label{lem_vecteur_G_inv}If $\tilde{Y}\in\mathcal{V}\left(\mathcal{M}\right)$
is a symplectic lift of a vector field $Y\in\mathcal{V}\left(\mathcal{B}\right)$,
then it is $\mathcal{G}$-invariant.
\end{lem}
We now state the announced decomposition of symplectic vector fields.
We stress the fact that this result still holds in the presence of
monodromy

\begin{prop}
\label{prop_decomposition_vectors}Any symplectic vector field $X\in\mathcal{V}\left(\mathcal{M}\right)$
decomposes in a unique way as \[
X=X_{1}+X_{2},\]
where 
\begin{itemize}
\item $X_{1}$ is a Hamiltonian vector field, $X_{1}=X_{A}$, with $\left\langle A\right\rangle =0$,
where $\left\langle A\right\rangle $ is the vertical average of the
Hamiltonian $A$.
\item $X_{2}$ is a symplectic lift of a vector field on $\mathcal{B}$.
Namely, it is the vertical average of $X$, i.e. $X_{2}=\left\langle X\right\rangle $.
\end{itemize}
\end{prop}
\begin{proof}
Let $\alpha=\omega\left(X,.\right)$ be the 1-form associated with
$X$, which is closed since $X$ is symplectic. Let $\alpha_{2}=\left\langle \alpha\right\rangle $
be the vertical average of $\alpha$ and let $\alpha_{1}=\alpha-\alpha_{2}$.
The $1$-forms $\alpha_{1}$ and $\alpha_{2}$ are closed since $d\left\langle \alpha\right\rangle =\left\langle d\alpha\right\rangle $.
Thus, the vector fields $X_{1}$ and $X_{2}$, associated with $\alpha_{1}$
and $\alpha_{2}$, are symplectic. On the other hand, one has $\left\langle \alpha_{1}\right\rangle =0$
and Lemma \ref{lem_forme_ferme_moyenne_null_exacte} then implies
that $X_{1}$ is Hamiltonian, $X_{1}=X_{A}$, with $\left\langle A\right\rangle =0$.
Finally, $\left\langle \alpha_{2}\right\rangle =\alpha_{2}$ implies
that $\left\langle X_{2}\right\rangle =X_{2}$. Now, any $\mathcal{G}$-invariant
vector field must be a lift of a vector field on $\mathcal{B}$, since
the toric action of $\mathcal{G}$ is vertical and transitive on each
fiber. This proves the second point of the proposition.

Moreover, the decomposition $X=X_{1}+X_{2}$ is the unique one of
this type. Indeed, suppose that there is a second decomposition $X=X_{1}^{'}+X_{2}^{'}$
with the same properties. Taking the vertical average of both expressions,
we obtain $\left\langle X_{1}+X_{2}\right\rangle =\left\langle X_{1}^{'}+X_{2}^{'}\right\rangle $
and thus $\left\langle X_{2}\right\rangle =\left\langle X_{2}^{'}\right\rangle $.
Now, by Lemma \ref{lem_vecteur_G_inv}, both $X_{2}$ and $X_{2}^{'}$
are $\mathcal{G}$-invariant. It follows that $X_{2}=X_{2}^{'}$ and
thus $X_{1}=X_{1}^{'}$.
\end{proof}

\section{Deformations of completely integrable systems}

\subsection{Regular deformations of completely integrable systems}

Let $\left(H_{0},\mathcal{M}\overset{\pi}{\rightarrow}\mathcal{B}\right)$
be a \emph{regular CI system} composed of a fibration in Lagrangian
tori $\mathcal{M}\overset{\pi}{\rightarrow}\mathcal{B}$ and a Hamiltonian
$H_{0}\in C^{\infty}\left(\mathcal{M}\right)$ constant along the
fibers. As we mentioned in the introduction, we will restrict ourselves
to regular deformations of $H_{0}$, i.e. smooth families of Hamiltonians
$H_{\varepsilon}$ which are CI and regular for each $\varepsilon$.
At this point, we would like to stress the fact that this does \emph{not}
imply that $H_{\varepsilon}$ is constant along the fibers of a family
of fibrations $\mathcal{M}\overset{\pi_{\varepsilon}}{\rightarrow}\mathcal{B}$
\emph{depending smoothly on $\varepsilon$}. Nevertheless, we conjecture
that is is true for the generic class of \emph{non-degenerate} Hamiltonians.
We will discuss the nondegeneracy conditions in Section \ref{sec_ND}
and we now restrict our study to the following class of deformations. 

\begin{defn}
\label{def_deformation_reg}Let $\left(H_{0},\mathcal{M}\overset{\pi}{\rightarrow}\mathcal{B}\right)$
be a regular CI system and $H_{\varepsilon}\in C^{\infty}\left(\mathcal{M}\right)$
a smooth family of Hamiltonians. We say that $H_{\varepsilon}$ is
a \emph{regular deformation} of $H_{0}$ if it has the form\[
H_{\varepsilon}=I_{\varepsilon}\circ\phi^{\varepsilon},\]
where $I_{\varepsilon}\in\pi^{*}\left(C^{\infty}\left(\mathcal{B}\right)\right)$
is a smooth family of functions with $I_{0}=H_{0}$ and $\phi^{\varepsilon}:\mathcal{M}\rightarrow\mathcal{M}$
is a smooth family of symplectomorphisms with $\phi^{0}=\mathbb{I}$.
\end{defn}
For our purposes, we will need to work from now on with time-dependent
vector fields since each smooth family of diffeomorphisms $\phi^{\varepsilon}$
with $\phi^{0}=\mathbb{I}$ is the flow at time $\varepsilon$ of
the time-dependent vector field $X_{\varepsilon}$ defined by \[
\frac{d\left(f\circ\phi^{\varepsilon}\left(m\right)\right)}{d\varepsilon}=X_{\varepsilon}\left(f\right)\circ\phi^{\varepsilon}\left(m\right)\]
for each smooth function $f\in C^{\infty}\left(\mathcal{M}\right)$
and each point $m\in\mathcal{M}$. We denote this flow by $\phi_{X_{\varepsilon}}^{\varepsilon}$.
In all the following, all the considered families $\phi^{\varepsilon}$
of diffeomorphisms will implicitly depend smoothly on $\varepsilon$
and satisfy $\phi^{0}=\mathbb{I}$. We refer e.g. to \cite{marsden_ratiu}
for a review of the properties of time-dependent vector fields.

\subsection{Non-degenerate CI systems\label{sec_ND}}

Nondegeneracy conditions are those used in K.A.M. theories, like for
example those introduced by Arnol'd, Kolmogorov, Bryuno or Rüssmann.
We refer to \cite{roy_5} for a review of different nondegeneracy
conditions together with their properties and we will focus on two
of them. But first of all, we need to define a few notions. 

Since the CI Hamiltonian $H_{0}$ is constant along the fibers which
are connected, it must be of the form $H_{0}=F_{0}\circ\pi$, with
$F_{0}\in C^{\infty}\left(\mathcal{B}\right)$. Denote by $\nabla$
the Duistermaat's affine connection which exists naturally on the
base space $\mathcal{B}$. For any subset $\mathcal{O}\subset\mathcal{B}$,
we denote by $\mathcal{V}_{\nabla}\left(\mathcal{O}\right)$ the space
of parallel vector fields. Since the holonomy of $\nabla$ may not
vanish, the space $\mathcal{V}_{\nabla}\left(\mathcal{O}\right)$
might be empty. Nevertheless, when $\mathcal{O}$ is simply connected,
this space is a $d$-dimensional vector space. All the nondegeneracy
conditions, including those presented here, are local : $F_{0}$ (or
$H_{0}$) is said to be \emph{non-degenerate} if is non-degenerate
at each $b\in\mathcal{B}$. Moreover, these conditions involve the
space of parallel vector fields, but the mentioned local character
means that one needs actually only the spaces $\mathcal{V}_{\nabla}\left(\mathcal{O}\right)$
for a neighborhood $\mathcal{O}\subset\mathcal{B}$ of each point
$b\in\mathcal{B}$. We will use a slight misuse of language and say
{}``for each $X\in\mathcal{V}_{\nabla}\left(\mathcal{B}\right)$''
instead of {}``for each $b\in\mathcal{B}$, each neighborhood $\mathcal{O}\subset\mathcal{B}$
of $b$ and each $X\in\mathcal{V}_{\nabla}\left(\mathcal{O}\right)$''.

For each $X\in\mathcal{V}_{\nabla}\left(\mathcal{B}\right)$, let
us define the function $\Omega_{X}\in C^{\infty}\left(\mathcal{B}\right)$
by $\Omega_{X}=dF_{0}\left(X\right)$ and the associated \emph{resonance
set} \[
\Sigma_{X}=\left\{ b\in\mathcal{B}\mid\Omega_{X}\left(b\right)=0\right\} .\]

\begin{defn}
\textbf{\label{def_russmann}}The function $F_{0}$ is \textbf{}\emph{Rüssmann
non-degenerate} if for each non-vanishing $X\in\mathcal{V}_{\nabla}\left(\mathcal{B}\right)$,
the resonant set $\Sigma_{X}$ has an empty interior.
\end{defn}
Among the nondegeneracy conditions used in the literature, Rüssmann's
Condition \cite{russmann} is the weakest one and has nevertheless
the following important consequence (see e.g. \cite{roy_5} for a
proof).

\begin{lem}
\label{lem_russmann_ND_fibration_unique}If $\left(H_{0},\mathcal{M}\overset{\pi}{\rightarrow}\mathcal{B}\right)$
is a {}``Rüssmann'' non-degenerate C.I system, then $\mathcal{M}\overset{\pi}{\rightarrow}\mathcal{B}$
is the unique fibration in Lagrangian tori such that $H_{0}$ is constant
along the fibers.
\end{lem}
This nondegeneracy condition is enough to insure the unicity of the
normal form of Theorem \ref{theo_deformation_hamiltonian} which will
be proved in Section \ref{sec_normal_form} but for the study of first
order deformations developed in Section \ref{sec_first_order_deformations},
we will need a stronger one, which is nevertheless weaker than Kolmogorov's
or Arnold's ones. 

\begin{defn}
\label{def_weak_ND}The function $F_{0}$ is \textbf{}\emph{weakly
non-degenerate} if for each non-vanishing $X\in\mathcal{V}_{\nabla}\left(\mathcal{B}\right)$
and each point $b\in\Sigma_{X}$, one has \[
d\left(\Omega_{X}\right)_{b}\neq0.\]
 
\end{defn}
This condition implies among other that the resonant sets $\Sigma_{X}$
are $1$-codimen\-sionnal submanifolds of $\mathcal{B}$.

\subsection{Normal form for regular deformations\label{sec_normal_form}}

The aim of this section is to prove Theorem \ref{theo_deformation_hamiltonian}
which insures that, by changing the function $I_{\varepsilon}$, one
may assume that $\phi^{\varepsilon}$ is a Hamiltonian flow. This
result is based on Proposition \ref{prop_decomposition_flow} which
states that any family of symplectomorphisms $\phi^{\varepsilon}$
can be written as the composition of a Hamiltonian flow with a family
of fiber-preserving symplectomorphisms. Let us first define precisely
these two notions.

\begin{defn}
\label{def_deformation_hamiltonian}A family of symplectomorphisms
$\phi^{\varepsilon}$ is called \emph{Hamiltonian} if its vector field
$X_{\varepsilon}$ is Hamiltonian, $X_{\varepsilon}=X_{A_{\varepsilon}}$,
with $A_{\varepsilon}\in C^{\infty}\left(\mathcal{M}\right)$ depending
smoothly on $\varepsilon$. 
\end{defn}

\begin{defn}
A family of diffeomorphisms $\phi^{\varepsilon}:\mathcal{M}\rightarrow\mathcal{M}$
is called \emph{fiber-preserving} if there exists a family of diffeomorphisms
on the base space $\varphi^{\varepsilon}:\mathcal{B}\rightarrow\mathcal{B}$
such that\[
\pi\circ\phi^{\varepsilon}=\varphi^{\varepsilon}\circ\pi.\]
We say that $\phi^{\varepsilon}$ is \emph{vertical} whenever $\varphi^{\varepsilon}=\mathbb{I}$
for all $\varepsilon$.
\end{defn}
Whenever a vector field on $\mathcal{M}$ is both symplectic and a
lift of a vector field on $\mathcal{B}$, then we have the following
property.

\begin{lem}
\label{lem_vecteur_symplec_et_lift}If $\tilde{Y}_{\varepsilon}\in\mathcal{V}\left(\mathcal{M}\right)$
is symplectic for each $\varepsilon$ and is a lift of a time-dependent
vector field $Y_{\varepsilon}\in\mathcal{V}\left(\mathcal{B}\right)$,
then it is $\mathcal{G}$-invariant and for each tensor field $T$
one has \[
\left\langle \left(\phi_{\tilde{Y}_{\varepsilon}}^{\varepsilon}\right)_{*}T\right\rangle =\left(\phi_{\tilde{Y}_{\varepsilon}}^{\varepsilon}\right)_{*}\left\langle T\right\rangle .\]

\end{lem}
\begin{proof}
Let denote by $\phi^{\varepsilon}=\phi_{\tilde{Y}_{\varepsilon}}^{\varepsilon}$
the flow of $\tilde{Y}_{\varepsilon}$. This flow is fiber-preserving
and thus verifies $\pi\circ\phi^{\varepsilon}=\varphi^{\varepsilon}\circ\pi$
with $\varphi^{\varepsilon}:\mathcal{B}\rightarrow\mathcal{B}$ a
family of diffeomorphisms. One can easily show that $\varphi^{\varepsilon}$
is actually the flow of $Y_{\varepsilon}$. 

First of all, for each vertical and parallel vector field $X\in\Gamma\left(\bigcup_{b\in\mathcal{B}}\mathcal{V}_{\nabla}\left(\mathcal{M}_{b}\right)\right)$,
one has $\phi_{*}^{\varepsilon}X\in\Gamma\left(\bigcup_{b\in\mathcal{B}}\mathcal{V}_{\nabla}\left(\mathcal{M}_{b}\right)\right)$.
Indeed, as mentioned in Section \ref{sec_period_bundle}, $\phi_{*}^{\varepsilon}X$
is vertical and parallel if and only if the $1$-form $\omega\left(\phi_{*}^{\varepsilon}X\right)$
is a pull-back. Now, one has $\omega\left(\phi_{*}^{\varepsilon}X\right)=\left(\left(\phi^{\varepsilon}\right)^{-1}\right)^{*}\left(\omega\left(X\right)\right)$
since $\phi^{\varepsilon}$ is symplectic for each $\varepsilon$.
On the other hand, $\omega\left(X\right)=\pi^{*}\beta$ with $\beta\in\Omega^{1}\left(\mathcal{B}\right)$,
since by hypothesis $X$ is vertical and parallel. Consequently, one
has \[
\omega\left(\phi_{*}^{\varepsilon}X\right)=\left(\left(\phi^{\varepsilon}\right)^{-1}\right)^{*}\pi^{*}\beta=\pi^{*}\left(\left(\varphi^{\varepsilon}\right)^{-1}\right)^{*}\beta.\]
This proves that $\omega\left(\phi_{*}^{\varepsilon}X\right)$ is
a pull-back and therefore $\phi_{*}^{\varepsilon}X$ is vertical and
parallel. 

If in addition $X\in\Gamma\left(\mathcal{O},\Lambda\right)$, with
$\mathcal{O}\subset\mathcal{B}$ a subset, i.e. $X$ is $1$-periodic
in $\pi^{-1}\left(\mathcal{O}\right)$, then so is $\phi_{*}^{\varepsilon}X$
in $\phi^{\varepsilon}\left(\pi^{-1}\left(\mathcal{O}\right)\right)$.
Now, the smooth bundle $\Lambda$ has discrete fibers and $\phi_{*}^{\varepsilon}X$
depends smoothly on $\varepsilon$. This implies that for all $\varepsilon$,
one has $\phi_{*}^{\varepsilon}X=\phi_{*}^{\varepsilon=0}X$ and thus
$\phi_{*}^{\varepsilon}X=X$. Then, the derivative with respect to
$\varepsilon$ shows that $\left[\tilde{Y},X\right]=0$, i.e. $\tilde{Y}$
is $\mathcal{G}$-invariant. By linearity, this is true as well for
all $X\in\Gamma\left(\bigcup_{b\in\mathcal{B}}\mathcal{V}_{\nabla}\left(\mathcal{M}_{b}\right)\right)$. 

Therefore, for each $X\in\Gamma\left(\bigcup_{b\in\mathcal{B}}\mathcal{V}_{\nabla}\left(\mathcal{M}_{b}\right)\right)$
and each $\varepsilon$, $\phi^{\varepsilon}$ commutes with the flow
$\phi_{X}^{t}$. This implies that $\phi^{\varepsilon}$ commutes
with the toric action of $\mathcal{G}$ and thus with the averaging
process, i.e. \[
\left\langle \left(\phi_{\tilde{Y}_{\varepsilon}}^{\varepsilon}\right)_{*}T\right\rangle =\left(\phi_{\tilde{Y}_{\varepsilon}}^{\varepsilon}\right)_{*}\left\langle T\right\rangle \]
for any tensor field $T$. 
\end{proof}
We can now give the following decomposition result for families of
symplectomorphisms.

\begin{prop}
\label{prop_decomposition_flow}Each family of symplectomorphisms
$\phi^{\varepsilon}$ decomposes in a unique way as follows : \[
\phi^{\varepsilon}=\Phi^{\varepsilon}\circ\phi_{Z_{\varepsilon}}^{\varepsilon},\]
 where 
\begin{itemize}
\item $\Phi^{\varepsilon}$ is a fiber-preserving family of  symplectomorphisms.
\item $Z_{\varepsilon}=X_{G_{\varepsilon}}$ is a time-dependent Hamiltonian
vector field with $\left\langle G_{\varepsilon}\right\rangle =0$.
\end{itemize}
Moreover, the vector field of $\Phi^{\varepsilon}$ is equal to the
average $\left\langle X_{\varepsilon}\right\rangle $, where $X_{\varepsilon}$
is the vector field of $\phi^{\varepsilon}$. 
\end{prop}
\begin{proof}
Let $X_{\varepsilon}$ be the vector field of $\phi^{\varepsilon}$.
Proposition \ref{prop_decomposition_vectors} insures that for each
$\varepsilon$, $X_{\varepsilon}$ decomposes into $X_{\varepsilon}=\tilde{Y}_{\varepsilon}+W_{\varepsilon}$,
where $\tilde{Y}_{\varepsilon}$ is a lift of a vector field $Y_{\varepsilon}\in\mathcal{V}\left(\mathcal{B}\right)$
and $W_{\varepsilon}$ is Hamiltonian. Moreover, by looking more carefully
at the proof of Proposition \ref{prop_decomposition_vectors}, one
can easily check that $\tilde{Y}_{\varepsilon}$ and $W_{\varepsilon}$
depend smoothly on $\varepsilon$, since $\tilde{Y}_{\varepsilon}$
is nothing but the vertical average of $X_{\varepsilon}$.

Let $\Psi^{\varepsilon}$ be the family of symplectomorphisms defined
by $\phi_{\tilde{Y}_{\varepsilon}+W_{\varepsilon}}^{\varepsilon}=\phi_{\tilde{Y}_{\varepsilon}}^{\varepsilon}\circ\Psi^{\varepsilon}$
and let $Z_{\varepsilon}$ be its vector field. On the one hand, $\Phi^{\varepsilon}=\phi_{\tilde{Y}_{\varepsilon}}^{\varepsilon}$
is fiber-preserving since $\tilde{Y}_{\varepsilon}$ is a lift of
a vector field on $\mathcal{B}$. On the other hand, one can check
in a straightforward way that the vector field $X_{\varepsilon}^{3}$
of a composition of flows $\phi_{X_{\varepsilon}^{1}}^{\varepsilon}\circ\phi_{X_{\varepsilon}^{2}}^{\varepsilon}$
is given by the formula $X_{\varepsilon}^{3}=X_{\varepsilon}^{1}+\left(\phi_{X_{\varepsilon}^{1}}^{\varepsilon}\right)_{*}X_{\varepsilon}^{2}$.
Therefore, in our case we have $\tilde{Y}_{\varepsilon}+W_{\varepsilon}=\tilde{Y}_{\varepsilon}+\phi_{\tilde{Y}_{\varepsilon}}^{\varepsilon}\left(Z_{\varepsilon}\right)$
and thus \[
Z_{\varepsilon}=\left(\phi_{\tilde{Y}_{\varepsilon}}^{\varepsilon}\right)_{*}^{-1}\left(W_{\varepsilon}\right).\]
According to Proposition \ref{prop_decomposition_vectors}, $W_{\varepsilon}$
is Hamiltonian and verifies $\left\langle W_{\varepsilon}\right\rangle =0$.
First, this insures that $Z_{\varepsilon}$ is Hamiltonian. Second,
Lemma \ref{lem_vecteur_symplec_et_lift} implies that \[
\left\langle Z_{\varepsilon}\right\rangle =\left(\phi_{\tilde{Y}_{\varepsilon}}^{\varepsilon}\right)_{*}^{-1}\left\langle W_{\varepsilon}\right\rangle =0\]
since $\tilde{Y}_{\varepsilon}$ is symplectic and a lift of a vector
field on $\mathcal{B}$.

Finally, we show that this decomposition is unique. Indeed, suppose
that we have a second decomposition $\phi_{X_{\varepsilon}}^{\varepsilon}=\phi_{\tilde{Y}_{\varepsilon}^{'}}^{\varepsilon}\circ\phi_{Z_{\varepsilon}^{'}}^{\varepsilon}$
with the same properties. The vector field $\tilde{Y}_{\varepsilon}^{'}\textrm{ }$
must be a lift of a vector field on $\mathcal{B}$ since $\phi_{\tilde{Y}_{\varepsilon}^{'}}^{\varepsilon}$
is fiber-preserving. On the other hand, as we mentionned before, we
have the relation $\tilde{X}_{\varepsilon}=\tilde{Y}_{\varepsilon}^{'}+\phi_{\tilde{Y}_{\varepsilon}^{'}}^{\varepsilon}\left(Z_{\varepsilon}^{'}\right)$.
Arguing as before, we can show that $\phi_{\tilde{Y}_{\varepsilon}^{'}}^{\varepsilon}\left(Z_{\varepsilon}^{'}\right)$
is a Hamiltonian vector field with vanishing vertical average. Now,
Theorem \ref{prop_decomposition_vectors} tells us that the decomposition
$X_{\varepsilon}=\tilde{Y}_{\varepsilon}+W_{\varepsilon}$ is unique
and thus $\tilde{Y}_{\varepsilon}^{'}=\tilde{Y}_{\varepsilon}$ and
$Z_{\varepsilon}^{'}=Z_{\varepsilon}$.
\end{proof}
We have now all the necessary material to state the following theorem
which gives a normal form for regular deformations of a given regular
CI system. 

\begin{thm}
\label{theo_deformation_hamiltonian}Let $\left(H_{0},\mathcal{M}\overset{\pi}{\rightarrow}\mathcal{B}\right)$
a regular CI system. If $H_{\varepsilon}$ is a regular deformation
of $H_{0}$, then there exist a family of functions $I_{\varepsilon}\in\pi^{*}\left(C^{\infty}\left(\mathcal{B}\right)\right)$
and a family of Hamiltonian symplectomorphisms $\phi_{X_{G_{\varepsilon}}}^{\varepsilon}$,
with $\left\langle G_{\varepsilon}\right\rangle =0$ such that\[
H_{\varepsilon}=I_{\varepsilon}\circ\phi_{X_{G_{\varepsilon}}}^{\varepsilon}\]
for each $\varepsilon$.

Moreover, if $H_{0}$ is Rüssmann non-degenerate, then the families
$I_{\varepsilon}$ and $\phi_{X_{G_{\varepsilon}}}^{\varepsilon}$
are unique.
\end{thm}
\begin{proof}
By definition, $H_{\varepsilon}$ is a regular deformation of $H_{0}$
if there exist a family of functions $J_{\varepsilon}\in\pi^{*}\left(C^{\infty}\left(\mathcal{B}\right)\right)$
and a family of symplectomorphisms $\phi^{\varepsilon}$ such that
$H_{\varepsilon}=J_{\varepsilon}\circ\phi^{\varepsilon}$. On the
other hand, Proposition \ref{prop_decomposition_flow} insures that
$\phi^{\varepsilon}$ decomposes into $\phi^{\varepsilon}=\Phi^{\varepsilon}\circ\phi_{X_{G_{\varepsilon}}}^{\varepsilon}$,
where $\Phi^{\varepsilon}$ is fiber-preserving and $\left\langle G_{\varepsilon}\right\rangle =0$.
Therefore, we have $H_{\varepsilon}=I_{\varepsilon}\circ\phi_{X_{G_{\varepsilon}}}^{\varepsilon}$,
where the function $I_{\varepsilon}=J_{\varepsilon}\circ\Phi^{\varepsilon}$
is indeed an element of $\pi^{*}\left(C^{\infty}\left(\mathcal{B}\right)\right)$
since $\Phi^{\varepsilon}$ is fiber-preserving.

Let us now show the unicity in case $H_{0}$ is Rüssmann non-degenerate.
Suppose there is another family of functions $I_{\varepsilon}^{'}\in\pi^{*}\left(C^{\infty}\left(\mathcal{B}\right)\right)$
and another family of symplectomorphims $\phi_{X_{G_{\varepsilon}^{'}}}^{\varepsilon}$,
with $\left\langle G_{\varepsilon}^{'}\right\rangle =0$ and such
that $H_{\varepsilon}=I_{\varepsilon}^{'}\circ\phi_{X_{G_{\varepsilon}^{'}}}^{\varepsilon}$.
We thus have $I_{\varepsilon}\circ\phi_{X_{G_{\varepsilon}}}^{\varepsilon}=I_{\varepsilon}^{'}\circ\phi_{X_{G_{\varepsilon}^{'}}}^{\varepsilon}$
and if we define the flow $\Phi^{\varepsilon}=\phi_{X_{G_{\varepsilon}}}^{\varepsilon}\circ\left(\phi_{X_{G_{\varepsilon}^{'}}}^{\varepsilon}\right)^{-1}$,
then we obtain $I_{\varepsilon}\circ\Phi^{\varepsilon}=I_{\varepsilon}^{'}$. 

First of all, since $\Phi^{\varepsilon}$ is a family of symplectomorphisms,
the fibration $\mathcal{M}\overset{\pi_{\varepsilon}}{\rightarrow}\mathcal{B}$
given by $\pi_{\varepsilon}=\pi\circ\left(\Phi^{\varepsilon}\right)^{-1}$
is also Lagrangian. Then, we can see that the function $I_{\varepsilon}$
is also constant along the fibers of the deformed fibration. Indeed,
by hypothesis it has the form $I_{\varepsilon}=I_{\varepsilon}^{'}\circ\left(\Phi^{\varepsilon}\right)^{-1}$.
Using then the fact that $I_{\varepsilon}^{'}$ has the form $I_{\varepsilon}^{'}=f_{\varepsilon}\circ\pi$
with $f_{\varepsilon}\in C^{\infty}\left(\mathcal{B}\right)$, it
follows that $I_{\varepsilon}=f_{\varepsilon}\circ\pi_{\varepsilon}$.
This proves that $I_{\varepsilon}$ is constant along the fibers of
both fibrations $\pi$ and $\pi_{\varepsilon}$.

Moreover, for $\varepsilon=0$ the function $I_{\varepsilon}$ is
equal to $H_{0}$ which is non-degenerate. This implies that $I_{\varepsilon}$
is also non-degenerate for small enough $\varepsilon$ since nondegeneracy
is an open condition. Therefore, Lemma \ref{lem_russmann_ND_fibration_unique}
insures that there is a unique fibration such that $I_{\varepsilon}$
is constant along the fibers. The two fibrations $\pi$ and $\pi_{\varepsilon}$
thus coincide, this proves that $\Phi^{\varepsilon}$ preserves the
initial fibration $\pi$. Consequently, we have the decomposition
$\phi_{X_{G_{\varepsilon}}}^{\varepsilon}=\Phi^{\varepsilon}\circ\phi_{X_{G_{\varepsilon}^{'}}}^{\varepsilon}$
with $\Phi^{\varepsilon}$ preserving the fibration $\pi$. Now, Proposition
\ref{prop_decomposition_flow} insure that this decomposition is unique.
Accordingly we have $G_{\varepsilon}=G_{\varepsilon}^{'}$ and thus
$I_{\varepsilon}=I_{\varepsilon}^{'}$.
\end{proof}

\subsection{First order deformations\label{sec_first_order_deformations}}

In this last section, we adress the problem of finding what are the
necessary and sufficient conditions on a perturbation $H_{1}\in C^{\infty}\left(\mathcal{M}\right)$
which insure that the perturbed Hamiltonian $H_{\varepsilon}=H_{0}+\varepsilon H_{1}$
is CI up to $\varepsilon^{2}$, i.e. has the form $H_{\varepsilon}=I_{\varepsilon}\circ\phi_{X_{G_{\varepsilon}}}^{\varepsilon}+O\left(\varepsilon^{2}\right)$
with $I_{\varepsilon}\in\pi^{*}\left(C^{\infty}\left(\mathcal{B}\right)\right)$
and $\left\langle G_{\varepsilon}\right\rangle =0$. 

Most of the work here will be achieved with the help of Fourier series.
Let us begin by expliciting the geometric status of the object we
will consider%
\footnote{We refer to \cite{roy_1} for a detailled description of this issue.%
}. First, we will work locally in some $\mathcal{O}\subset\mathcal{B}$,
with an action-angle coordinates system $\left(\xi,x\right)$ and
consider the Fourier series with respect to the periodic variable
$x$. For any smooth function $f\left(\xi,x\right)$ we will denote
by $\tilde{f}\left(\xi,k\right)$ its Fourier series defined by the
usual expression \[
f\left(\xi,x\right)=\sum_{k\in E}e^{ik\left(x-x_{0}\right)}\tilde{f}\left(\xi,k\right).\]
The discrete set $E$ in which the Fourier variable $k$ lives is
naturally a lattice of the vector space $\mathcal{V}_{\nabla}\left(\mathcal{O}\right)$
of parallel vector fields on $\mathcal{O}$. This can be seen as follows.
First, if $\xi$ denotes the coordinates of a point $b\in\mathcal{O}$,
then $x-x_{0}$ can be understood as an element of $\mathcal{V}_{\nabla}\left(\mathcal{M}_{b}\right)$
well-defined up to elements of $\Lambda_{b}$, the lattice of $1$-periodic
parallel vector fields on the fiber $\mathcal{M}_{b}$. On the other
hand, its dual $\Lambda_{b}^{*}$ is a lattice of the space $\Omega_{\nabla}^{1}\left(\mathcal{M}_{b}\right)$
of parallel $1$-forms on $\mathcal{M}_{b}$. The Fourier variable
$k$ lives naturally in $\Lambda_{b}^{*}$. Moreover, this family
$\Lambda_{b}^{*}$ depends smoothly on $b$ as $\Lambda_{b}$ does.
Now, the symplectic form provides an isomorphism $\iota_{b}:\Omega_{\nabla}^{1}\left(\mathcal{M}_{b}\right)\rightarrow T_{b}\mathcal{B}$,
depending smoothly on $b$. Under this identification, $k$ can be
seen as an element of the lattice $\iota_{b}\left(\Lambda_{b}^{*}\right)$
of the vector space $T_{b}\mathcal{B}$. Finally, if we identify $\iota_{b}\left(\Lambda_{b}^{*}\right)$
with the space $E\subset\mathcal{V}_{\nabla}\left(\mathcal{O}\right)$
of sections of the associated lattice bundle $\bigcup_{b\in\mathcal{O}}\iota_{b}\left(\Lambda_{b}^{*}\right)$,
we obtain a suitable space for the Fourier variable $k$ to live in.
Accordingly, for each $k\in E$ the Fourier series $\tilde{f}\left(\xi,k\right)$
is a smooth (with respect to $\xi$) function, well-defined up to
a phase, due to an arbitrary choice of the family of origin points
$b\rightarrow x_{0}\left(b\right)$. 

\begin{defn}
\label{def_non_resonant_function}A function $f\in C^{\infty}\left(\mathcal{M}\right)$
is called \emph{non-resonant} if for each non-vanishing $k\in E$
and each $\xi\in\Sigma_{k}$ one has $\tilde{f}\left(\xi,k\right)=0$. 
\end{defn}
The resonant manifolds $\Sigma_{k}$ were defined in Section \ref{sec_ND}.
We have the following equivalent criterion which has to be checked
on each torus $\mathcal{M}_{b}$ on which the dynamics of $X_{H_{0}}$
is periodic.

\begin{lem}
\label{lem_non_resonant_function}A function $f\in C^{\infty}\left(\mathcal{M}\right)$
is non-resonant \emph{}if and only if for each $T$-periodic torus
$\mathcal{M}_{b}$ the average of $f$ along the trajectories of $X_{H_{0}}$
\[
\overline{f}:=\frac{1}{T}\int_{0}^{T}\left.f\right|_{\mathcal{M}_{b}}\circ\phi_{X_{H_{0}}}^{t}dt\]
is a constant function on $\mathcal{M}_{b}$.
\end{lem}
\begin{proof}
First, one can show (see e.g. \cite[Prop. A.62]{roy_1}) that for
each $k\neq0$ the set of periodic tori in $\Sigma_{k}$ is dense
in $\Sigma_{k}$. This implies that the nonresonance condition is
equivalent to \[
\forall\xi\textrm{ periodic},\forall k\in E\setminus0,dF_{0}\left(k\right)_{\xi}=0\Longrightarrow\widetilde{f}\left(\xi,k\right)=0,\]
where $F_{0}\in C^{\infty}\left(\mathcal{B}\right)$ is the function
defined by $H_{0}=F_{0}\circ\pi$. On the other hand, a short calculation
shows that for each $\xi$ the Fourier series $\tilde{\bar{f}}\left(\xi,k\right)$
of the average $\bar{f}$ is given by \[
\widetilde{\overline{f}}\left(\xi,k\right)=\left\{ \begin{array}{c}
\widetilde{f}\left(x,k\right)\textrm{ if }dF_{0}\left(k\right)=0\textrm{ }\\
0\textrm{ if }dF_{0}\left(k\right)\neq0\end{array}\right..\]
Therefore, the nonresonance condition indeed amounts to requiring
that the averaged function $\bar{f}$ is constant on the torus $\mathcal{M}_{b}$.

This nonresonance condition is the right one to controls the complete
integrability up to $\varepsilon^{2}$, as it is shown in the following
theorem.
\end{proof}
\begin{thm}
\label{theo_first_order_deformation}Let $\left(H_{0},\mathcal{M}\overset{\pi}{\rightarrow}\mathcal{B}\right)$
be a weakly non-degenerate regular CI system and $H_{1}\in C^{\infty}\left(\mathcal{M}\right)$
a perturbation. The perturbed Hamiltonian $H_{\varepsilon}=H_{0}+\varepsilon H_{1}$
is CI up to $\varepsilon^{2}$ if and only if $H_{1}$ is non-resonant.
\end{thm}
\begin{proof}
First, the complete integrability up to $\varepsilon^{2}$ means that
$H_{\varepsilon}$ has the form $H_{\varepsilon}=I_{\varepsilon}\circ\phi_{X_{G_{\varepsilon}}}^{\varepsilon}+O\left(\varepsilon^{2}\right)$
with $I_{\varepsilon}\in\pi^{*}\left(C^{\infty}\left(\mathcal{B}\right)\right)$
and $\left\langle G_{\varepsilon}\right\rangle =0$. In this expression,
the terms of order $\varepsilon^{0}$ give simply $H_{0}=I_{0}$ and
the $\varepsilon^{1}$ terms yield the equation $H_{1}=I_{1}+X_{G_{0}}\left(H_{0}\right)$.
By definition of the Poisson bracket, this is equivalent to \begin{equation}
\left\{ H_{0},G_{0}\right\} =I_{1}-H_{1}.\label{equ_1}\end{equation}
In Fourier coordinates, this equation reads \[
idF_{0}\left(k\right)\tilde{G}_{0}\left(\xi,k\right)=\tilde{I}_{1}\left(\xi,k\right)-\tilde{H}_{1}\left(\xi,k\right)\]
for each $\xi$ and each $k\in E$. For $k=0$, we have $\tilde{I}_{1}\left(\xi,0\right)=I_{1}\left(\xi\right)$
since $I_{1}$ is a function constant along the fibers. We thus have
to set $\tilde{H}_{1}\left(\xi,0\right)=I_{1}\left(\xi\right)$. Remark
that $\tilde{H}_{1}\left(\xi,0\right)$ is nothing but the vertical
average of the function $H_{1}$. The Fourier coefficient $\tilde{G}_{0}\left(\xi,0\right)$
is free and can be set to $0$, which means that $\left\langle G_{0}\right\rangle =0$.
Now, for all non-vanishing $k$, we have $\tilde{I}_{1}\left(\xi,k\right)=0$
and we need to solve the equation \[
dF_{0}\left(k\right)\tilde{G}_{0}\left(\xi,k\right)=i\tilde{H}_{1}\left(\xi,k\right).\]
 The nonresonance condition is certainly necessary, since in order
to divide by the function $\Omega_{k}=dF_{0}\left(k\right)$, $\tilde{H}_{1}\left(\xi,k\right)$
needs to vanish at least where $\Omega_{k}$ does, i.e. on the resonance
manifold $\Sigma_{k}$. The solution $\tilde{G}_{0}$ is thus defined
by the quotient $\frac{\tilde{H}_{1}\left(\xi,k\right)}{\Omega_{k}\left(\xi\right)}$
and it still remains to prove that the nonresonance condition is sufficient
to insure that $\tilde{G}_{0}\left(\xi,k\right)$ is smooth with respect
to $\xi$ uniformly with respect to $k$, and with a fast decay in
$k$. This will mean that $G_{0}\left(\xi,x\right)$ is smooth with
respect to $\left(\xi,x\right)$. For this purpose, we will show that
for any compact set $\mathcal{K}\subset\mathcal{B}$, there are two
positive constants $T$ and $C$ such that \begin{equation}
\left|\Omega_{k}\right|<T\Longrightarrow\left|d\Omega_{k}\right|>C\label{equ_2}\end{equation}
uniformly with respect to $k\in E\setminus0$. If this holds, then
for each $k$ we decompose $\mathcal{K}$ into two parts defined by
$\left|\Omega_{k}\right|<T$ and $\left|\Omega_{k}\right|\geq T$.
Away from the resonance manifold $\Sigma_{k}$, i.e for $\xi$ such
that $\left|\Omega_{k}\left(\xi\right)\right|\geq T$, we can simply
devide by $\Omega_{k}$ and $\tilde{G}_{0}$ will satisfy the estimate
$\left|\tilde{G}_{0}\left(\xi,k\right)\right|\leq\frac{\tilde{H}_{1}\left(\xi,k\right)}{T}$.
On the other hand, close to the resonance manifold $\Sigma_{k}$,
i.e for $\xi$ such that $\left|\Omega_{k}\left(\xi\right)\right|<T$,
we consider $X_{k}$ the gradient of $\Omega_{k}$ (for some fixed
riemannain metric). Its norm verifies $\left|X_{k}\right|>C$ and
is transversal to the submanifold $\Sigma_{k}$. It is thus suitable
to parametrize the {}``distance'' to $\Sigma_{k}$, thru its flow
$\phi_{X_{k}}^{t}$. Indeed, since $X_{k}$ is the dual vector of
$\Omega_{k}$, i.e. $d\Omega_{k}\left(X_{k}\right)=1$, we have $\Omega_{k}\circ\phi_{X_{k}}^{t}\left(\xi\right)=t+\Omega_{k}\left(\xi\right)$.
Therefore we can compute the first order Taylor expansion with integral
rest of $\tilde{H}_{1}\left(\xi,k\right)$. This yields \[
\tilde{H}_{1}\left(\xi,k\right)=\tilde{H}_{1}\left(\phi_{X_{k}}^{-\Omega_{k}\left(\xi\right)},k\right)+\int_{0}^{\Omega_{k}\left(\xi\right)}dt\, X_{k}\left(\tilde{H}_{1}\right)\circ\phi_{X_{k}}^{t-\Omega_{k}\left(\xi\right)}\left(\xi\right).\]
By construction, the point $\phi_{X_{k}}^{-\Omega_{k}\left(\xi\right)}$
is on $\Sigma_{k}$ and $\tilde{H}_{1}$ vanishes at this point according
to the nonresonance condition. Then, the change of variable $u=t/\Omega_{k}\left(\xi\right)$
gives \[
\tilde{H}_{1}\left(\xi,k\right)=\Omega_{k}\left(\xi\right)\int_{0}^{1}du\, X_{k}\left(\tilde{H}_{1}\right)\circ\phi_{X_{k}}^{\left(u-1\right)\Omega_{k}\left(\xi\right)}\left(\xi\right)\]
 and we can solve the equation $\Omega_{k}\left(\xi\right)\tilde{G}_{0}\left(\xi,k\right)=i\tilde{H}_{1}\left(\xi,k\right)$
by dividing by $\Omega_{k}$. The fast decay of $\tilde{H}_{1}\left(\xi,k\right)$
implies the fast decay of the solution $\tilde{G}_{0}\left(\xi,k\right)$
for both cases $\left|\Omega_{k}\right|<T$ and $\left|\Omega_{k}\right|\geq T$,
and thus proves the smoothness of $G_{0}\left(\xi,x\right)$.

The last point is to prove the existence of the constants $T$ and
$C$ in Equation \ref{equ_2}. In fact we will prove that this equation
holds for $k$ living in the space $P=\left\{ X\in\mathcal{V}_{\nabla}\left(\mathcal{B}\right),\left|X\right|\geq1\right\} $
and this will imply the result for $k\in E\setminus0$. Because of
nondegeneracy, one has $d\Omega_{k}\neq0$ and thus $\frac{d\Omega_{k}}{\left|k\right|}\neq0$
for all point in $\Sigma_{k}\cap\mathcal{K}$. Therefore, there is
a constant $C\left(k\right)$ such that $\left|d\Omega_{\frac{k}{\left|k\right|}}\right|=\frac{\left|d\Omega_{k}\right|}{\left|k\right|}>2C\left(k\right)$
in $\Sigma_{k}\cap\mathcal{K}$. Now, the smoothness of $\Omega_{\frac{k}{\left|k\right|}}$
implies that there is a constant $T\left(k\right)$ such that $\frac{\left|d\Omega_{k}\right|}{\left|k\right|}>C\left(k\right)$
whenever $\left|\Omega_{\frac{k}{\left|k\right|}}\right|<T\left(k\right)$.
Let us now decompose the elements $k$ into their angular and radial
parts, i.e. $k:=\left(\frac{k}{\left|k\right|},k\right)\in S^{d-1}\times\left[1,\infty\right]$.
Taking the minimum of $T\left(k\right)$ and $C\left(k\right)$ over
the compact set $S^{d-1}$, we obtain positive constants $T^{'}\left(\left|k\right|\right)$
and $C^{'}\left(\left|k\right|\right)$ such that $\frac{\left|d\Omega_{k}\right|}{\left|k\right|}>C^{'}\left(\left|k\right|\right)$
whenever $\left|\Omega_{\frac{k}{\left|k\right|}}\right|<T^{'}\left(\left|k\right|\right)$.
Using again $\Omega_{\frac{k}{\left|k\right|}}=\frac{\Omega_{k}}{\left|k\right|}$
and setting $T=T^{'}\left(1\right)$ and $C=C^{'}\left(1\right)$,
we see that the following implication holds.

\[
\left|\frac{\Omega_{k}}{\left|k\right|}\right|<T\Longrightarrow\frac{\left|d\Omega_{k}\right|}{\left|k\right|}>C.\]
Finally, using the fact that $\left|k\right|\leq1$, we obtain Equation
\ref{equ_2}.
\end{proof}

\end{document}